\documentclass[12pt,a4paper]{amsart}
\usepackage{enumerate}
\usepackage{geometry}
\usepackage{amsmath}
\usepackage{amssymb}
\usepackage{amsthm}
\usepackage{amsrefs}
\usepackage{color}
\usepackage{xspace}
\usepackage{comment}
\usepackage[driverfallback=dvipdfm]{hyperref}
\newtheorem{theorem}{Theorem}
\newtheorem{lemma}[theorem]{Lemma}

\def\blem{\begin{lemma}}
\def\elem{\end{lemma}}
\def\bmat{\begin{pmatrix}}
\def\emat{\end{pmatrix}}

\def\ep{\varepsilon}
\def\R{\mathbb{R}}
\def\N{\mathbb{N}}
\def\ol{\overline}
\def\fr{\frac}
\def\mid{\,:\,} 
\def\gl{\lambda}
\def\IN{\text{ in }}

\def\AND{\text{ and }}
\def\FOR{\text{ for }}

\def\ON{\text{ on }}
\def\bproof{\begin{proof}}
\def\eproof{\end{proof}}
\def\stm{\setminus}

\def\pl{\partial}
\def\gd{\delta}
\def\tim{\times}
\def\bald{\begin{aligned}}
\def\eald{\end{aligned}} 
\def\FORALL{\text{ for all }}

\DeclareMathOperator\USC{USC}
\DeclareMathOperator\LSC{LSC}

\DeclareMathOperator\tr{tr}

\def\stm{\setminus}
\def\cS{\mathcal{S}}

\def\fr{\frac} 
  
\def\ga{\alpha}     
\def\go{\omega}
      
\def\ep{\gep}    
\def\mid{\,:\,}   
\def\gb{\beta} 

\def\gd{\delta}
\def\gz{\zeta} 
\def\gth{\theta}   
 
\def\gl{\lambda}
\def\gL{\Lambda}
\def\gs{\sigma}   
                  
\def\tim{\times}

\def\ol{\overline}
\def\ul{\underline}           
\def\pl{\partial}

\def\gG{\varGamma}

\def\bcases{\begin{cases}}
\def\ecases{\end{cases}}
\def\balns{\begin{align*}}
\def\ealns{\end{align*}}
\def\bald{\begin{aligned}}
\def\eald{\end{aligned}}
\def\beq{\begin{equation}}
\def\eeq{\end{equation}}
\def\bred{\begin{color}{red}} \def\ered{\end{color}}
\def\gO{\Omega} 

\def\1{\mathbf{1}}

\def\IN{\text{ in }} \def\FOR{\text{ for }} 
\def\AND{\text{ and }}

\def\ON{\text{ on }}

\theoremstyle{definition}
\newtheorem{definition}{Definition}%[section]
\theoremstyle{plain}
\newtheorem{proposition}[definition]{Proposition}
\theoremstyle{remark}
\newtheorem{remark}[definition]{Remark}

\def\dirichlet{(D$_\ep$)\xspace}
\def\dirichletzero{(D$_0$)\xspace}
\def\ga{\gamma}

\allowdisplaybreaks[4]

\def\dirichlet{(D$_\ep$)\xspace}

\def\dirichletzero{(D$_0$)\xspace}

\title[Thin domains with oblique boundary condition]{Fully nonlinear elliptic PDEs in thin domains with oblique-Dirichlet mixed boundary conditions}

\author[I. Birindelli]{Isabeau Birindelli}
\address[I. Birindelli]{Dipartimento di Matematica Guido Castelnuovo, Sapienza 
Universit\`a di Roma, Piazzale Aldo Moro 5, Roma, Italy.}
\email{isabeau@mat.uniroma1.it}

\author[A. Briani]{Ariela Briani}
\address[A. Briani]{Université de Tours, 
Université d’Orléans, CNRS, IDP, UMR 7013, Tours, France.}
\email{ariela.briani@univ-tours.fr}

\author[H. Ishii]{Hitoshi Ishii%$^*$
}
\address[H. Ishii]{Institute for Mathematics and Computer Science, Tsuda  University,
 2-1-1 Tsuda, Kodaira, Tokyo 187-8577 Japan.}
\email{hitoshi.ishii@waseda.jp}
\keywords{asymptotic behaviour of solutions, thin domains}

\thanks{A. Briani was partially supported by l’Agence Nationale de la Recherche (ANR), project
ANR-22-CE40-0010 COSS; I. Birindelli was partially supported by project Leoni 2023 GNAMPA-INDAM and project  "At the Edge of Reaction-diffusion equations" Sapienza Università di Roma. H. {Ishii was partially supported by the JSPS KAKENHI Grant No. JP25K07072.}} %The project was very much advanced while I. Birindelli was visiting prof. Ishii in Tsuda University, she wished to thank the Institution for the invitation.}

\subjclass[2020]{
35B40, %Asymptotic behavior of solutions,
35D40, %Viscosity solutions 
35J25  	%Boundary value problems for second-order elliptic equations
49L25 %Viscosity solutions
}

\def\gO{\Omega}\def\ep{\varepsilon}\def\T{\mathrm{T}}\def\B{\mathrm{B}}
\def\L{\mathrm{L}} \def\sT{\intercal}
\def\cM{\mathcal{M}}
\def\ga{\alpha}
\def\Tilde{\widetilde}
\def\Hat{\widehat}
\begin{document}

\begin{abstract} 
We consider asymptotic behavior of solutions to  the oblique-Dirichlet mixed boundary conditions without the \emph{strict monotonicity} of the equation in the variable corresponding to the unknown function for "thin domains" i.e.  when the $N+1$ dimensional domains collapse to an $N$ dimensional domain. A global ellipticity condition in the limit equation is introduced.
\end{abstract}

\maketitle 

\tableofcontents
\allowdisplaybreaks

\section{Introduction} 
In this paper, we investigate the asymptotic behavior of solutions to fully nonlinear elliptic partial differential equations (PDEs) in thin domains with oblique and Dirichlet boundary conditions. Specifically, we consider the following boundary value problem \dirichlet, with $0<\ep\ll 1$, on a domain $\Omega_\epsilon \subset \Omega \times [-1, 1]\subset \R^{N+1}$

\beq \tag{D$_\ep$} \label{dirichlet}
\left\{\ 
\bald
&F(D^2u^\ep,Du^\ep,u^\ep,z)=0 \ \IN\gO_\ep, 
\\
&\gamma^+\cdot Du^\ep=\beta^+ \ \ \ON\pl_\T\gO_\ep,  
\\
&\gamma^-\cdot Du^\ep=\beta^- \ \ \ON \pl_\B\gO_\ep,
\\
&u^\ep(z)=\beta(z) \ON\pl_\L\gO_\ep,
\eald
\right.
\eeq 
where
 the domain 
$\gO_\ep$ is defined as
\[
\gO_\ep=\{z=(x,y)\in \gO\tim \R\mid \ep g^-(x)<y<\ep g^+(x)\},
\] 
 $\gO$ is a given bounded domain of $\R^N$ and $g^\pm$, 
 are given continuous functions on $\ol\gO$
satisfying $g^-<g^+$ on $\ol\gO$. Here, 
the boundary $\pl\gO_\ep$ is divided into three parts, called the top, bottom, and lateral boundaries, which are denoted by $\pl_\T\gO_\ep$, $\pl_\B\gO$, 
and $\pl_\L\gO$, respectively, and they are defined as 
\[\bald
&\pl_\T\gO=\{(x,y)\mid x\in\ol\gO,\,y=\ep g^+(x)\},
\quad
\\&
\pl_\B\gO_\ep=\{(x,y)\mid x\in\ol\gO,\, y=\ep g^-(x)\},
\quad
\\&\pl_\L\gO_\ep=\{(x,y)\mid x\in\pl\gO,\, \ep g^-(x)\leq y\leq \ep g^+(x)\}. 
\eald
\]
Furthermore, 
$F$ is a Bellman-Isaacs type operator, on $ \cS(N+1)\tim\R^{N+1}\tim\R\tim \ol\gO\tim[-1,1]$,  
%\iffalse %%%%%
 of the form 
\beq\label{eq5.5+}
F(X,p,r,z)=\inf_{\gl \in L}\sup_{\mu\in M}(-\tr \gs_{\gl\mu}^\sT\gs_{\gl\mu} X-b_{\gl\mu}\cdot p+c_{\gl\mu} r-f_{\gl\mu}),
\eeq
where $\cS(N+1)$ denotes the space of symmetric matrices of size $N+1$, 
$L$ and $M$ are given nonempty sets,
and $\gs_{\gl\mu}$, 
 $b_{\gl\mu}$, $c_{\gl\mu}$, and $f_{\gl\mu}$ are functions on 
$\ol{\gO}\tim[-1.1]$ 
taking values in $k\tim(N+1)$ matrices for some $k\in\N$, 
$\R^{N+1}$, $\R$, and $\R$, respectively. 
 %%\[
%%F_{\gl\mu}(X,p,r,z)=-\tr \gs_{\gl\mu}^\sT\gs_{\gl\mu} X-b_{\gl\mu}\cdot p+c_{\gl\mu} r-f_{\gl\mu},
%%\]
%with 
%\[ \bald
%&\gs_{\gl\mu}\in C(\ol\gO\tim[-1,1],\cM(k,N+1)),\quad 
%b_{\gl\mu}\in C(\ol\gO\tim[-1,1],\R^{N+1}), 
%\quad
%\\& c_{\gl\mu},\,f_{\gl\mu}\in C(\ol\gO\tim[-1,1],\R).
%\eald
%\]
%In the above, 
%$\cM(m,n)\subset\R^{m\tim n}$ denotes the set of all real $m\tim n$ matrices,
%and $\cS(m)$ the set of all symmetric matrices in $\cM(m,m)$,  
%and $k\in\N$ is a fixed number throughout this paper.   {\color{magenta}  It is not clear who is $k$ %here? Is defined below?}

Further 
%\fi %%%%
and 
detailed requirements on $F$ and $\gamma^\pm, \beta^\pm, \beta$ will be provided later on 
(see \eqref{h.1}, \eqref{h.2}, \eqref{non-negative}, \eqref{cont}, 
\eqref{h.4}, \eqref{h.5}, and \eqref{h.6} below). 

This work is a direct continuation of the Dirichlet lateral boundary condition part of our previous study (see \cite{BBI2}), where we established results for similar problems under the assumption of strict monotonicity of the operator $F(X,p,r,z)$ in the variable $r$. In the present paper, we relax this requirement. Instead of strict monotonicity, we assume only that $F$ is non-decreasing in $r$, which corresponds to the condition $c_{\lambda\mu}(z) \geq 0$ for all $(\lambda, \mu) \in L \times M$ (see \eqref{non-negative}).

Let us mention that typical condition for existence of solutions for Neumann boundary problem associated with fully nonlinear degenerate operators requires strict monotonicity in the 
dependence on the solution see e.g. \cite{BNe,INe}. Here, the lateral condition is a Dirichlet condition and this gives some room. Nonetheless to compensate for the lack of strict 
monotonicity, we introduce a requirement on the ellipticity of $F$ in the limit domain $\Omega$. This request which reinforce the degenerate ellipticity is much less than the strict ellipticity.
%\iffalse %%%%

The precise condition is that there exists a function $s \in C^2(\Omega, \mathbb{R})$ such that for all $x \in \Omega$,
\begin{equation} \label{cond_9}
\inf_{(\gl,\mu)\in L\tim M,\ x\in\ol\gO}\, \left|
(Ds(x), -Ds(x)\gamma_0^T(x)) A_{\lambda\mu}(x, 0)\right|> 0,
\end{equation}
where $A_{\lambda\mu} := \sigma_{\lambda\mu}^\sT \sigma_{\lambda\mu}$ and 
$\gamma_0\in C(\ol\gO,\R^N)$ .
%\fi %%%%%
%which will be formulated %defined 
%shortly (see \eqref{positive} below). 
\begin{remark}\label{rem-A-positive} The existence of $s\in C^2(\ol{\gO},\R)$ in \eqref{cond_9}
or \eqref{G-positive} is a global condition. To simplify the situation, we consider a smooth function 
$A: \ol\gO\to\cS(N)$ satisfying $A(x)\geq 0$. If 
\beq\label{A-positive}
A(x)\not=0\quad\FOR \ x\in\ol\gO, 
\eeq
then, for each point $y\in\ol\gO$, one can choose a linear function $s_y$ for which 
$Ds_y(x)A(x)Ds_y^\sT(x)>0$ in a neighborhood of $y$. But, in general, condition \eqref{A-positive}
does not imply the existence of a $C^2$-function $s$ on $\ol\gO$ such that 
\beq\label{A-positive*}
Ds(x)A(x)Ds^\sT(x)>0\quad\FOR x\in\ol\gO . 
\eeq
\end{remark}
%%%{\color{magenta} T}
The existence of such a function $s$ over the entire domain $\Omega$ is non-trivial, more is explained in Section \ref{example} to illustrate this phenomenon. 
%An example showing this phenomenon is given in Section \ref{example}.  
%
%As it will be discussed in Remark \ref{rem-A-positive} and Section \ref{example}, this condition has a non-local character. While a local version might follow from the non-zero nature of $A_{\gl\mu}(x)$, the existence of such a function $s$ over the entire domain $\Omega$ is non-trivial, and we provide a counterexample in Section \ref{example} to illustrate this phenomenon. 

Under these revised hypotheses, our main result, Proposition \ref{main.prop}, establishes the existence and uniform boundedness of viscosity solutions $u^\ep$ for sufficiently small $\ep > 0$. Our proof relies on the construction of appropriate strict sub- and super-solutions (Lemma \ref{main.lem}).  In the proof of Lemma \ref{main.lem}, we employ a coordinate transformation to distort the system so that the boundary conditions match the simplified case where $\gamma_0 = 0$. Once the existence result of Proposition \ref{main.prop} is proved, to obtain our main convergence result for the limiting thin domain, we only need to follow the convergence argument in \cite{BBI1,BBI2}.

\section{Standing assumptions} \label{hypo}
In this section, we collect the assumptions on the problem \dirichlet, most of them are taken from
\cite{BBI2}.

$\bullet$ Assumptions on $F$:\  $F$ has the form: 
\beq\label{h.1} %\beq  \label{eq5.5+}
F(X,p,r,z)=\inf_{\gl \in L}\sup_{\mu\in M} F_{\gl\mu}(X,p,r,z)
\ \ \ON \cS(N+1)\tim\R^{N+1}\tim\R\tim \ol\gO\tim[-1,1],  
\eeq
where 
$L,\,M$ are given nonempty sets and 
\[
F_{\gl\mu}(X,p,r,z)=-\tr \gs_{\gl\mu}^\sT\gs_{\gl\mu} X-b_{\gl\mu}\cdot p+c_{\gl\mu} r-f_{\gl\mu},
\]
with 
\[ \bald
&\gs_{\gl\mu}\in C(\ol\gO\tim[-1,1],\cM(k,N+1)),\quad 
b_{\gl\mu}\in C(\ol\gO\tim[-1,1],\R^{N+1}), 
\quad
\\& c_{\gl\mu},\,f_{\gl\mu}\in C(\ol\gO\tim[-1,1],\R).
\eald
\]
In the above, $\cM(m,n)\subset\R^{m\tim n}$ denotes the set of all real $m\tim n$ matrices,
and $\cS(m)$ the set of all symmetric matrices in $\cM(m,m)$,  
and $k\in\N$ is a fixed number throughout this paper.   
The following boundedness 
and continuity are required: 
\beq\label{h.2}\left\{\,\bald
&\max\{|\gs_{\gl\mu}(z)|,|b_{\gl\mu}(z)|,|c_{\gl\mu}(z)|,|f_{\gl\mu}(z)|\}\leq C_F,
\\&\max\{|\gs_{\gl\mu}(z)-\gs_{\gl\mu}(z')|,|b_{\gl\mu}(z)-b_{\gl\mu}(z')|\}\leq C_F|z-z'|,
\\&\max\{|c_{\gl\mu}(z)-c_{\gl\mu}(z')|,|f_{\gl\mu}(z)-f_{\gl\mu}(z')|\}\leq \go_F(|z-z'|), 
\eald\right.
\eeq
where $C_F$ is a positive constant, $\go_F$ is a modulus of continuity. 
Additionally, 
\beq\label{non-negative}
c_{\gl\mu}(z)\geq 0 \ \ \FOR z\in\ol\gO\tim[-1,1]. 
\eeq

$\bullet$ Assumptions on $\gO$: 
The set $\gO\subset\R^N$ is a bounded $C^1$ domain and satisfies the 
exterior sphere condition, i.e., 
%%%\begin{minipage}{0.85\textwidth}%% $\pl\gO$ is $C^1$ and 
there is a constant $r_0>0$ such that
\beq\label{h.3} %\left\{
B_{r_0}(x+r_0\nu(x))\cap\gO=\emptyset \ \ \ \FOR x\in\pl\gO. 
\eeq
%%\end{minipage}%%%\right.
%%\]
where $\nu(x)$ denotes the unit outward normal at $x$.

$\bullet$ Assumptions on $\gb,\,\gamma^\pm,\,\beta^\pm$: \ The following 
continuity condition is assumed.
\beq\label{cont}\left\{\bald
&\gb\in C(\ol\gO\tim[-1,1],\R),\qquad
\gamma^\pm=(\gamma_1^\pm,\gamma_2^\pm)\in C(\ol\gO\tim[-1,
%%%%{\color{magenta}, }
1],\R^{N+1}),
\quad 
\\& \beta^\pm \in C(\ol\gO\tim[-1,1],\R). 
\eald\right. 
\eeq
It is natural to assume that $\pm\gamma_2^\pm >0$ on $\ol\gO\tim[-1,1]$, and in fact, by normalization, we assume that 
\beq\label{h.4}
\gamma_2^\pm=\pm 1. 
\eeq
Furthermore, we assume that there exist functions 
\beq \label{h.5} 
\gamma_o\in C^2(\ol\gO,\R^N), \quad\beta_o\in C^2(\ol\gO,\R),
\quad 
k^\pm\in C^1(\ol\gO,\R^N),\ \ \AND \ \  l^\pm\in C(\ol\gO,\R),
\eeq
such that 
\beq\label{h.6}
\left\{\,\bald
&\gamma_1^\pm(x,y)=\pm\gamma_o(x)+k^\pm(x)y+o(|y|),
\\&\beta^\pm(x,y)=\pm\beta_o(x)+l^\pm(x)y+o(|y|),
\eald \right.
\eeq
as $y\to 0$, where $o(|y|)/|y|\to 0$ uniformly on $\ol\gO$ as $y\to 0$.
These imply:   
\beq\label{h.7} 
\beta^+(x,0)=-\beta^-(x,0)=:\beta_o(x)\ \ \AND \ \ 
\gamma_1^+(x,0)=-\gamma_1^-(x,0)=:\gamma_o(x) \ \FOR x\in\ol\gO. 
\eeq
In \cite{BBI2}, the regularity assumption on $\gamma_0,\,\beta_0$ is that 
$(\gamma_0,\beta_0)\in C^1(\ol\gO,\R^{N+1})$. 
Here, we assume \eqref{h.5}, a stronger regularity assumption 
on $\gamma_0,\,\beta_0,\, k^\pm$ 
than that in \cite{BBI2},  for brevity of presentation.

$\bullet$ Ellipticity assumptions: We set $A_{\gl\mu}:=\gs_{\gl\mu}^\sT\gs_{\gl\mu}$. 
The assumptions are: 
\beq  \label{positive}\left\{\,\begin{minipage}{0.85\textwidth}
there 
exists a function $s\in C^2(\ol\gO,\R)$ such that  for $x\in\ol\gO$, 
\[
\inf_{(\gl,\mu)\in L\tim M}
\left|(Ds(x),-Ds(x)\gamma_0^\sT(x))A_{\gl\mu}(x,0)\right|>0.
\]
\end{minipage}\right.
\eeq
Furthermore for any $x\in \pl\gO$, 
\beq\label{b-positive}
\inf_{(\gl,\mu)\in L\tim M}\left|(\nu(x), -\nu(x)\gamma_0^\T(x)) A_{\gl,\mu}(x,0)\right|>0.  
\eeq
This condition and \cite[(73)]{BBI2} appears to be different, but are actually the same.
Indeed, recall here the condition \cite[(73)]{BBI2}, which  is stated as follows: for any $x\in\pl\gO$,
\beq\label{b-positive*}
\inf_{(\gl,\mu)\in L\tim M}|\nu(x)\bmat I_N & -\gamma_0^\sT(x)\emat \gs_{\gl\mu}^\sT(x,0)|>0.
\eeq
We set $
v:=\nu(x)\bmat I_N&-\gamma_0^\sT(x)\emat=(\nu(x),-\nu(x)\gamma_0^\sT(x))\in\R^{N+1}, 
$
and note that 
\[
|v \gs_{\gl\mu}^\sT(x,0)|^2=v\gs_{\gl\mu}^\sT(x,0)\left(v\gs_{\gl\mu}^\sT(x,0)\right)^\sT
=v A_{\gl\mu}(x,0)v^\sT\leq |vA_{\gl\mu}(x,0)||v|,
\]
which shows that \eqref{b-positive*} implies \eqref{b-positive}. 
On the other hand, since 
\[|vA_{\gl\mu}(x,0)|\leq|v\gs_{\gl\mu}^\sT(x,0)|
|\gs_{\gl\mu}(x,0)|,\] 
we find that \eqref{b-positive} and \eqref{h.2} yields \eqref{b-positive*}. 

$\bullet$ Assumptions on $g^\pm$: \
We assume:
\beq\label{h.8}    
g^\pm\in C^1(\ol\gO),\quad g^-(x)<g^+(x) \ \ \FOR x\in\ol\gO.  
\eeq

\section{Around the ellipticity condition} \label{example}
Before entering into the proofs and statements of the core results of this paper we wish to dwell on the ellipticity condition in order to 
show how the condition \eqref{positive} affects the choice of the matrices $A_{\lambda\mu}(x,0)\geq 0$ allowed.

So precisely we consider $\cS(N)$-valued function $A$ on
$\ol\gO$, having the property, $A(x)\geq 0$ for $x\in\ol\gO$, such that there exists a function $s\in C^1(\ol\gO,\R)$
\beq \label{280326} Ds(x)A(x)Ds(x)^T> 0 . \eeq

It is quite obvious that if $A(x)$ is any non negative matrix such that at least one eigenvalue is positive, and its corresponding eigenvector $v(x)$ 
is a gradient field in $\Omega$ with potential $s(x)$ then,  \eqref{280326} is satisfied.

On the other hand having one positive eigenvalue is not enough. As mentioned in Remark \ref{rem-A-positive}, 
we now give an example of an $\cS(N)$-valued function $A$ on
$\ol\gO$, with one eigenvalue equal to $1$ for any $x\in\ol\gO$, but such that
 \eqref{280326} is not satisfied for any $s\in C^1(\ol\gO,\R)$.

Let $\gO$ be a domain in $\R^2$ that includes the unit circle $\pl B_1(0)$.  
For $\gth\in\R$, set 
\[
\bald
R(\gth)&:=\bmat\cos\gth&-\sin\gth \\ \sin\gth&\cos\gth\emat, 
\qquad
\\
A(\gth)&:= R(\gth)\bmat0&0\\0&1\emat R(\gth)^{-1}= R(\gth)\bmat0&0\\0&1\emat R(\gth)^{\sT}. 
\eald
\]
For the moment, we consider $A$ as an $\cS(2)$-valued function on $\pl B_1(0)$ by associating 
$A(\gth)$ to point $(\cos\gth,\sin\gth)\in\pl B_1(0)$.
Set
\[
v(\gth):=(1, 0) R(\gth)^{-1}=(\cos\gth,\sin\gth),
\qquad w(\gth):=(0,1)R(\gth)^{-1}=(-\sin\gth,\cos\gth),
\]
and note that $v(\gth), w(\gth)$ are (unit) eigenvectors of $A(\gth)$. Indeed, 
$v(\gth)A(\gth)=0$ and $w(\gth)A(\gth)=w(\gth)$, and the associated eigenvalues are $0$ and $1$, respectively. 
It is clear that $A(\gth)$ are all non-negative definite and non-zero (i.e., \eqref{A-positive} is satisfied). 

To proceed, we pick any $s\in C^1(\pl B_1(0),\R)$. Let $x(\gth)=(\cos \gth,\sin\gth)$, and 
note that
\[
\fr{d}{d\gth}s(x(\gth))=Ds(x(\gth))\cdot \dot x(\gth)
=Ds(x(\gth))\cdot (-\sin\gth,\cos\gth))
=Ds(x(\gth))\cdot w(\gth). 
\]
Now, since $s(x(0))=s(x(\pi))$, by the Mean Value Theorem, there is a $\gth_0\in (0,2\pi)$ such that  $(d/d\gth) s(x(\gth))\big|_{\gth=\gth_0}=0$. Thus, we have 
\[
Ds(x(\gth_0))\cdot w(\gth_0)=0, 
\]
which implies, together with the two-dimensionality of the space, that 
$Ds(x(\gth_0))=rv(\gth_0)$ for some $r\in\R$ and moreover, $Ds(x(\gth_0))A(\gth_0)=0$.  
This shows that for the $\cS(2)$-valued function $A$ on $\pl B_1(0)$, \eqref{A-positive*} does not hold for any  $s\in C^1(\pl B_1(0),\R)$.  It is worth to remark that the function $A$ can be extended 
smoothly and boundedly to $\R^2$ so that $A(x)\in\cS(2)$ and $A(x)\geq 0$ for $x\in\R^2$.  
To see this, let $\chi\in C^\infty(\R,\R)$ be such that 
$\chi(r)\geq 0$ for $r\in\R$, $\chi(r)=1$ for $r<1/2$, and $\chi(r)=0$ for $r\geq 1$. Then, for 
$x=r(\cos\gth,\sin\gth)\in\R^2$, define
\[
A(x)=R(\gth)\bmat \chi(r) & 0\\ 0&1\emat R(\gth)^{-1}.
\]  
It is clear this $A$ is an extension of the original $A$ defined on $\pl B_1(0)$.

Noting that $A(x)=\bmat0&0\\0&1\emat $ for $x\in B_{1/2}(0)$,
we easily deduce that $A$ is well-defined and smooth in $\R^2$. 
It is clear that $A(x)\in\cS(2)$ for $x\in\R^2$ and $\chi(|x|)$ and $1$ are the eigenvalues of $A(x)$, which implies that $A$ is bounded in $\R^2$ and non-negative definite everywhere. 

We wish to mention that the example was inspired by \cite{FS}.

\section{Main result}
\subsection{Limit problem} \label{Limit-prob}

Although the settings are slightly different from those in \cite{BBI2}, it will turn out that the limit problem, i.e., the 
PDE problem for the limit function of $u^\ep$ to \dirichlet, is identical.   
Recall (see \cite[(1.13), (1.14)]{BBI2}) that 
the limit problem is stated as
\begin{equation}\tag{D$_0$}
\left\{\,\bald
&G(D^2u,Du,u,x) =0 \quad \IN\gO, 
\\&u=\beta(x,0) \quad \ON\pl\gO.
\eald
\right.
\end{equation}
Here, $G$ has the form:  for 
$(X,p,r,x)\in\cS(N)\tim\R^N\tim\R\tim\ol\gO$,
\beq\label{l.1}\left\{\, \bald
&G(X,p,r,x):=\inf_{\gl\in L}\sup_{\mu\in M}G_{\gl\mu}(X,p,r,x),
\\
&G_{\gl\mu}(X,p,r,x):=-\tr (\Tilde A_{\gl\mu}(x)X)-\tilde b_{\gl\mu}(x)\cdot p
+\tilde c_{\gl\mu}(x)r-\tilde f_{\gl\mu}(x).
\eald\right.
\eeq
The coefficients in the above are as follows: It has been shown (see \cite{BBI2}) that for all 
$(X,p,r,x)\in \cS(N)\tim\R^N\tim\R\tim\ol\gO$,
\[
G(X,p,r,x)=F(A+B+C,(p,\gb_0-\gamma_0\cdot p),r,(x,0)),
\]
where 
\[\bald
A&=\bmat X&-X\gamma_0^\sT(x)\\ -\gamma_0(x)X&\gamma_0(x) X \gamma_0^\sT(x)\emat,
\\
B&=\bmat 0&-(p D\gamma_0(x))^\sT\\ -pD\gamma_0(x)&b(x)\cdot p\emat, 
%%%% {\color{magenta} , }
\\
C&=\bmat0&D\gb_0^\sT(x)\\ D\gb_0(x) &c(x)\emat,
\\
b&=\gamma_0 D\gamma_0^\sT -\fr{1}{g^+-g^-}(g^+k^++g^-k^-),
\\
c&=-\gamma_0 D\gb_0+\fr{1}{g^+-g^-}(g^+l^++g^-l^-).
\eald\]
(The functions $\gamma_0,\,\beta_0,\,k^\pm,\, l^\pm$ are those in \eqref{h.6}.) 
It follows from the above that 
\beq\label{l.2}\left\{\bald
\Tilde A_{\gl\mu}(x)&= \bmat I_N&-\gamma_0^\sT(x)\emat A_{\gl\mu}(x,0)\bmat I_N&-\gamma_0^\sT(x)\emat^\sT,
\\
\tilde b_{\gl\mu}(x)&=b_{\gl\mu}(x,0)\bmat I_N&-\gamma_0^\sT(x)\emat^\sT
\\&\quad -2e_{N+1}  A_{\gl\mu}(x,0)\bmat D\gamma_0(x)&0\emat^\sT
+e_{N+1} A_{\gl\mu}(x,0)e_{N+1}^\sT b(x),
\\
\tilde c_{\gl\mu}(x)&=c_{\gl\mu}(x,0),
\\
\tilde f_{\gl\mu}(x)&=f_{\gl\mu}(x,0) 
+b_{\gl\mu}(x,0)\gb_0(x)+\tr A_{\gl\mu}(x,0)\bmat 0& D\gb_0^\sT(x) \\
D\gb_0(x)&c(x)\emat.
\eald\right.
\eeq 
Here,  $e_{N+1}$ denotes the unit vector in $\R^{N+1}$, with its $(N+1)^{\,\text{th}}$ entry 
being unity. Setting 
\[
\tilde\gs_{\gl\mu}(x)=\gs_{\gl\mu}(x,0)\bmat I_N &-\gamma_0^\sT(x)\emat^\sT \ \ \FOR x\in\ol\gO,
\]
we have $\Tilde A_{\gl\mu}=\tilde\gs_{\gl\mu}^\sT\tilde\gs_{\gl\mu}$ on $\ol\gO$.
It follows from \eqref{h.2}, \eqref{h.5}, and \eqref{h.8}  
that for $x,x'\in\ol\gO$ and $(\gl,\mu)\in L\tim M$,
\beq\label{l.3}\left\{\,\bald
&\max\{|\tilde\gs_{\gl\mu}(x)|,|\tilde b_{\gl\mu}(x)|,|\tilde c_{\gl\mu}(x)|,|\tilde f_{\gl\mu}(z)|\}\leq \Tilde C_F,
\\&\max\{|\tilde\gs_{\gl\mu}(x)-\tilde\gs_{\gl\mu}(x')|,|\tilde b_{\gl\mu}(x)-\tilde b_{\gl\mu}(x')|\}\leq \Tilde C_F|x-x'|,
\\&\max\{|\tilde c_{\gl\mu}(x)-\tilde c_{\gl\mu}(x')|,|\tilde f_{\gl\mu}(x)-\tilde f_{\gl\mu}(x')|\}\leq \tilde \go_F(|x-x'|), 
\eald\right.
\eeq
where $\Tilde C_F$ is a positive constant and $\tilde \go_F$ is a modulus of continuity.

\subsection{Convergence result}
Our main result in this paper is as follows: 

\begin{theorem} \label{main.thm} 
Assume the hypotheses stated above. Then,  
\begin{enumerate}[(i)]
\item  there exists $\ep_1\in (0,1)$ such that for any $\ep\in(0,\ep_1)$, 
\dirichlet has a bounded viscosity solution on $\ol{\gO_\ep}$. 
\item There exists a unique viscosity 
solution $u^0\in C(\ol\gO)$ of \dirichletzero that satisfies the Dirichlet boundary 
condition in the classical sense (i.e., the pointwise sense). 
\item If $S_\ep$ be the collection of all bounded viscosity solutions to \dirichlet, then  
\[
\lim_{\ep\to 0^+}\sup_{u\in S_\ep} \sup_{(x,y)\in\ol{\gO_\ep}}|u(x,y)-u^0(x)|=0.
\]
\end{enumerate}
\end{theorem}

\section{A key proposition for the oblique-Dirichlet problem} \label{key-prop}
In our study of \dirichlet, 
unlike in the previous work \cite{BBI1,BBI2}, we do not assume strict monotonicity of $F(X,p,r,z)$ in $r$, 
and, by assuming instead \eqref{positive}, which is a sort of ellipticity of $F$,    
we establish a proposition similar to \cite[Proposition 2]{BBI2}.

Since we are interested in the asymptotic behaviour, as $\ep\to 0^+$, of the solutions to the boundary value problem \eqref{dirichlet},  we may focus our 
considerations on small $\ep$, say $0<\ep<\ep_0$, which guarantees that  
\beq \label{k.1}
\ol{\gO_\ep}\subset \ol\gO\tim[-1,1]. 
\eeq

We remark that the inequality in \eqref{positive} is equivalent that for $x\in\ol\gO$,
\beq\label{positive*}
\inf_{(\gl,\mu)\in L\tim M}(Ds(x),-Ds(x)\gamma_0^\sT(x))A_{\gl\mu}(x,0)(Ds(x),-Ds(x)\gamma_0^\sT(x))^\sT>0.
\eeq 
As we have seen the equivalence of \eqref{b-positive} and \eqref{b-positive*}, 
we can easily check that \eqref{positive} is equivalent to requiring that for $x\in\ol\gO$,
\[
\inf_{(\gl,\mu)\in L\tim M}\left|(Ds(x),-Ds(x)\gamma_0^\sT(x))\gs_{\gl\mu}^\sT(x,0)\right|>0.
\]
Writing $v:=(Ds(x),-Ds(x)\gamma_0^\sT(x))$, we observe that 
\[
|v\gs_{\gl\mu}^\sT(x,0)|^2=v\gs_{\gl\mu}^\sT(x,0)\left(v\gs_{\gl\mu}^\sT(x,0)\right)^\sT
%%%5{\color{magenta} \sT ???}
 =v A_{\gl\mu}(x,0)v^\sT.
\]
It is now obvious that \eqref{positive} and \eqref{positive*} are equivalent. 
Noting that the function: 
\[
\ol\gO\ni x\mapsto \inf_{(\gl,\mu)\in L\tim M}(Ds(x),-Ds(x)\gamma_0^\sT(x))A_{\gl\mu}(x,0)(Ds(x),-Ds(x)\gamma_0^\sT(x))^\sT
\]
is continuous, we find that \eqref{positive} implies
\beq \label{all-inf}
\inf_{x\in\ol\gO,\ (\gl,\mu)\in L\tim M}(Ds(x),-Ds(x)\gamma_0^\sT(x))A_{\gl\mu}(x,0)(Ds(x),-Ds(x)\gamma_0^\sT(x))^\sT>0.
\eeq

A simple condition stronger than \ref{positive} is  
that there is a constant, non-negative, and non-zero matrix $A_N\in\cS(N)$ such that  
for any $(\gl,\mu)\in L\tim M$,
\beq\label{k.2}
A_{\gl\mu}\geq 
\bmat
A_N &0\\
0&0
\emat. 
\eeq
Indeed, under the hypothesis 
%%%{\color{magenta} 
\eqref{k.2}, 
%%}, 
let $\xi\in\R^N\stm\{0\}$ be 
an eigenvector of $A_N$ associated with a non-zero eigenvalue. With the linear function $s(x):=\xi\cdot x$, the condition \eqref{positive} is satisfied. 

It is worth writing the condition \eqref{positive*} from the viewpoint of the limit problem \dirichletzero,   
where $G$ is the function describing the PDE in $\gO$, 
$G$ is ``superposed'' through $\inf$-$\sup$ operation as a family 
of affine functions $G_{\gl\mu}$, and the second-order term of the differential 
operator $G_{\gl\mu}$ has the coefficient $\Tilde A_{\gl\mu}$ (see \eqref{l.1}) 
Furthermore, as in \eqref{l.2}, we have
\[
\Tilde A_{\gl\mu}(x)=\bmat I_N&-\gamma_0^\sT(x)\emat A_{\gl\mu}(x,0)
\bmat I_N& -\gamma_0^\sT(x)\emat^\sT. 
\] 
Thus, the condition \eqref{positive*} can be stated as
\beq\label{G-positive}
\inf_{(\gl,\mu)\in L\tim M}Ds(x)\Tilde A_{\gl\mu}(x)Ds^\sT(x)>0 \ \ \FOR x\in\ol\gO.
\eeq
This is an ellipticity requirement for the equation $G=0$ .

%%%%%%%%%%%%%%%%%%
\iffalse
Furthermore, we assume that there exists a function $s:\ol\gO\to\R$ such that 
\beq\label{eq5.6+2}
s\in C^2(\ol\gO),\qquad \gamma_0\cdot Ds=0 \ \ \ON \ol\gO,\qquad |Ds|\geq 1 \ \ \ON\ol\gO. 
\eeq
This is a crucial reuirement in our approach. If $N\geq 2$ and $\gamma_0$ is a constant vecor, we choose a vector $v_0\in\R^N$ 
such that $v_0\cdot \gamma_0=0$ and $|v_0|=1$, set $s(x)=v_0\cdot x$ and then $s$ has the required properties \eqref{eq5.6+2}.  More generally, if $\gamma_0$ is a small and smooth 
perturbation of a constant vector, then we have a good chance to build a right $s(x)$ by 
solving the first order linear PDE $\gamma_0\cdot Ds=0$ via the characteristic method. 
By adding a constant, we may and do assume that $s\geq 0$ in $\gO$. If $\gO$ is the interval 
$(0,1)$, then the existence of such a function $s$ is possible only in the case when $\gamma_0=0$. 
If $N=2$ and $\gamma_0(x_1,x_2)=(x_2,-x_1)$, the flow generated by the vector field
$\gamma_0$ is the rotation around the origin and there is no smooth $s$ such that 
$\gamma_0\cdot Ds=0$ near the origin and $Ds(0)\not=0$.  
\fi
%%%%%%%%%%%%%%%%%%%%%%

The following is the main observation in this article, which is a version of \cite[Proposition 2]{BBI2}. 

\begin{proposition} \label{main.prop}
Under the above hypotheses, 
there exist positive constants $\ep_1<\ep_0$ and $C_0$ such that for each $0<\ep<\ep_1$, there is a viscosity solution  $u^\ep$ to  \eqref{dirichlet}. 
Furthermore any viscosity solution $v^\ep$ to \eqref{dirichlet} satisfies  
$\sup_{\ol{\gO_\ep}}|v^\ep|\leq C_0$. 
\end{proposition}

Our proof of the above proposition relies upon the following lemma. 

\begin{lemma}\label{main.lem} Under the same hypotheses of Proposition \ref{main.prop}, 
there exist constants $0<r<1$, $0<\ep_1<\ep_0$, $C>0$, and 
functions $\ol{\psi_\ep}, \ul{\psi_\ep}\in C^2(\ol{\gO}\tim[-r,r])$ such that 
for any $\ep\in(0,\ep_1)$, 
%%$\ol{\gO_\ep}\subset\ol\gO\tim[-r,r]$ and 
\begin{align}
\label{m.1}&\gamma^+\cdot D\ol{\psi_\ep}-\gb^+>0 \ \ \ON \pl_\T\gO_\ep,
\\\label{m.2}&\gamma^-\cdot D\ol{\psi_\ep}-\gb^->0 \ \ \ON \pl_\B\gO_\ep,
\\\label{m.3}&F(D^2\ol{\psi_\ep},D\ol{\psi_\ep},\ol{\psi_\ep},(x,y))>0 \ \ \ON \ol{\gO_\ep},
\\\label{m.4}&\gamma^+\cdot D\ul{\psi_\ep}-\gb^+<0 \ \ \ON \pl_\T\gO_\ep,
\\\label{m.5}& \gamma^-\cdot D\ul{\psi_\ep}-\gb^-<0 \ \ \ON \pl_\B \gO_\ep,
\\\label{m.6}&F(D^2\ul{\psi_\ep},D\ul{\psi_\ep},\ul{\psi_\ep},(z,y))<0 \ \ \ON \ol{\gO_\ep},
\\\label{m.7}&C>\ol{\psi_\ep}>\ul{\psi_\ep}>-C \ \ \ \ON \ol{\gO_\ep}. 
\end{align}
\end{lemma}

We accept this lemma and prove Proposition \ref{main.prop}. Then 
we give the proof of Lemma \ref{main.lem}. 

\bproof[Proof of Proposition \ref{main.prop}] Let $0<r<1$, $\ep_1$, $C>0$, and $\ol{\psi_\ep}, \ul{\psi_\ep}$, with 
$0<\ep<\ep_1$, be the constants and functions from Lemma \ref{main.lem}. 
By adding a positive constant to $\ol{\psi_\ep}$ and $-\ul{\psi_\ep}$ if necessary, we may 
assume that for $\ep\in(0,\ep_1)$, 
\beq\label{m.8}
\ol{\psi_\ep}>\gb>\ul{\psi_\ep} \ \ \ON \pl\gO\tim[-r,r]. 
\eeq
By 
%%the 
Perron's method, for each $\ep\in(0, \ep_1)$, we can construct a (viscosity) solution $u^\ep$ to \dirichlet, which satisfies $\ol{\psi}_\ep\leq u^\ep\leq \ol{\psi}_\ep$ on $\ol{\gO_\ep}$. 
%%%%{\color{magenta}.}

Let $\ep\in(0,\ep_1)$.  Let $v^\ep$ be a (viscosity) solution to \dirichlet. We claim that $\ul{\psi}_\ep\leq v^\ep\leq \ol{\psi}_\ep$ 
on $\ol{\gO_\ep}$. For the contradiction proof, 
suppose that for some $z\in\ol{\gO_\ep}$,
\[
\max((v^\ep)^*-\ol{\psi}_\ep)=((v^\ep)^*-\ol{\psi}_\ep)(z)>0,
\]
where $(v_\ep)^*$ denotes the upper semicontinuous envelope of $v^\ep$. 
Observe that if $z\in\gO_\ep$, then, by the subsolution property of $v^\ep$, 
we have
\beq\label{m.9}
0\geq F(D^2\ol{\psi_\ep}(z),D\ol{\psi_\ep}(z), (v^\ep)^*(z),z)
\geq F(D^2\ol{\psi_\ep}(z),D\ol{\psi_\ep}(z), \ol{\psi_\ep}(z),z),
\eeq
contradicting \eqref{m.3}.  
If $z\in \pl_\T\gO_\ep\stm\pl_\L\gO_\ep$, we have either \eqref{m.9} 
%\[0\geq F(D^2\ol{\psi_\ep}(z),D\ol{\psi_\ep}(z), (v^\ep)^*(z),z)
%\geq F(D^2\ol{\psi_\ep}(z),D\ol{\psi_\ep}(z), \ol{\psi_\ep}(z),z),
%\]
or
\beq\label{m.10}
0\geq \gamma^+(z)\cdot D\ol{\psi_\ep}(z)-\beta^+(z).
\eeq
But, these contradict \eqref{m.1} and \eqref{m.3}. 
The case when $z\in\pl_\B\gO_\ep$ can be treated similarly and leads a contradiction. 
If $z\in\pl_\L\gO_\ep\stm \left(\pl_\T\gO_\ep\cup\pl_\B\gO_\ep\right)$, then we have either \eqref{m.9}
or 
\[
0\geq (v^\ep)^*(z)-\beta(z)>\ol{\psi_\ep}(z)-\beta(z),
\]
leading a contradiction. Similarly, we can show that $v^\ep\geq \ul{\psi}_\ep$ on $\ol{\gO_\ep}$, 
and we 
%{\color{magenta} 
conclude
%} 
that $|v^\ep|\leq C$ on $\ol{\gO_\ep}$. 
\eproof 

\section{Proof of Lemma \ref{main.lem}} This section is fully devoted to the proof of Lemma \ref{main.lem}, which is divided into three steps.  

\subsection*{\quad Step 1} In this first step, we discuss the case where 
\beq\label{gamma00}
\gamma_0=0 \ \ \  \ON\ \ol{\gO}. 
\eeq
We need to build classical strict sub- and super-solutions to \eqref{dirichlet}. Choose a function  $h\in C^2(\ol\gO,\R)$ such that 
\[ g^-(x)<h(x)<g^+(x) \ \  \FOR x\in\ol\gO. 
\]
Thanks to \eqref{positive} and \eqref{positive*}, we have
\[
(Ds(x),0)A_{\gl\mu}(x,0)(Ds(x),0)^\sT>0 \ \ \FOR x\in\ol{\gO}. 
\]
By multiplying $s$ by a positive constant if necessary, we may assume that 
for some constant $0<r<1$ and 
for all $(x,y)\in\ol\gO\tim[-r,r]$,
\beq\label{positive**}
(Ds(x),0)A_{\gl\mu}(x,y)(Ds(x),0)^\sT\geq 1.
\eeq
Moreover, we may assume, by adding a constant to $s$, 
that $s\geq 0$ on $\ol\gO$. 

Let $\ga\geq 1$ be a constant to be fixed later, with the function $s$ of condition \eqref{positive} , we define function $\chi$ 
on $\ol\gO$ by $\chi(x)=e^{\ga s(x)}$. The choice of $\ga$ will be crucial in the following argument. 
Let $C_D$ be a positive constant to be fixed later. Note that
$1\leq \chi \leq e^{\ga\|s\|_\infty}=:C_\ga $, where $\|s\|_\infty=\max_{x\in\ol\gO}|s(x)|$. Set 
\[
\rho(x)=C_D+C_\ga -\chi(x),
\] 
and note that $C_D\leq \rho<C_D+C_\ga$ on $\gO$.

Let $\gL>1$ be a constant to be fixed later on. We define 
\[
\ol{\psi_\ep}(x,y)=\rho(x)+ \beta_0(x)y+\ga\gL\chi(x)(y-\ep h(x))^2,
\]
and
\[
\ul{\psi_\ep}(x,y)=-\rho(x)+\beta_0(x)y-\ga\gL\chi(x)(y-\ep h(x))^2.
\]
We will show that if $\ga$ and $\gL$ are large enough, then there exists $\ep_1\in(0,\ep_0)$ such that 
for any $0<\ep<\ep_1$, $\ol{\psi_\ep}$ and $\ul{\psi_\ep}$ are \emph{strict} 
super- and sub-solutions to \dirichlet, respectively.   

We temporarily introduce the following notational convention: for function $q=q(z,r)$ and positive function $p=p(r)$, we write 
\[
q(z,r)=O(p(r)) 
\] 
if there is a constant $C>0$, independent of $(z,r)$, such that 
\[
|q(z,r)|\leq Cp(r)
\]
for all $z$ and  $r$. Also, we write 
\[
q(z,r)=O^+(p(r))
\] 
if there are positive constants $C_0, C_1$, independent of $(z,r)$, such that 
\[
C_0\leq \fr{q(z,r)}{p(r)}\leq C_1
\]
for all $(z,r)$. The variable $z$ will be either $(x,y)$ or $x$, and $r$ will be either
$\ep$, $(\ep,\ga)$, or $(\ep,\ga,\gL)$.   

We will be only concerned with $\ol\psi_\ep$, and we write $\psi=\psi_\ep$ for 
$\ol\psi_\ep$. The function $\ul{\psi_\ep}$ can be treated in a parallel way. 
Regarding the boundary condition, we will be only concerned with the boundary 
$\pl_\T\gO\cup\pl_\L\gO_\ep$; the other part of the boundary can be treated similarly. 

In what follows we always assume that 
\beq\label{eq9}
0<\ep<1,\qquad \ga>1,\qquad \gL>1,\qquad \ep \ga<1,\qquad \ep\gL<1. 
\eeq

A comment on the choice of these constants may be in order. We are going to choose
the product $\ga\gL$ large enough for $\psi$ to satisfy the boundary condition on
$\pl_\T\gO_\ep$, which has an effect to let $\psi$ be convex at least in the $y$ variable, 
while we want $\psi$ to be somehow ``concave'' for $\psi$ being a strict supersolution in $\gO_\ep$.

We compute that 
\[\bald
D\rho&=-D\chi=-\ga \chi Ds, %%%\qquad Ds= (1,\ldots,1),
\\ D^2\chi&=\ga^2 \chi Ds\otimes Ds+\ga\chi D^2s=\ga^2 \chi Ds^\sT Ds+\ga\chi D^2s,
\eald
\]
which implies that on $\ol{\gO_\ep}$,
\[
D\rho=-D\chi=O(\ga)\chi,\qquad D^2\chi=O(\ga^2)\chi. 
\]

Now, we deal with the boundary condition on $\pl_\T\gO_\ep$. For this, let $(x,y)$ be any point in $\pl_\T\gO_\ep$, 
and note that $y=\ep g^+(x)=O(\ep)$, $\gamma_1^+=\gamma_0+O(\ep)=O(\ep)$ by \eqref{gamma00}, and $\gb^+=\gb_0+O(\ep)$. 
We compute that 
\beq\label{psi_x}
D_x\psi=-\ga \chi Ds +D\beta_0 y +\ga \gL D\chi (y-\ep h)^2-2\ep\ga\gL\chi(y-\ep h)Dh, 
\eeq
and 
\beq\label{psi_y}
D_y\psi=\beta_0+2\ga\gL\chi (y-\ep h).
\eeq
Noting that, since $y=\ep g^+(x)$ and $g^+(x)-h(x)\geq \min_{\ol\gO}(g^+-h)>0$, 
%%%on $\ol\gO$,  
\beq\label{O^+}
y-\ep h(x)=\ep(g^+(x)-h(x))=O^+(\ep),
\eeq
we deduce, in view of \eqref{eq9}, that
\[\bald
D_x\psi&=O(\ga)\chi+O(\ep)+O(\ga\gL \ep^2)|D\chi| +O(\ga\gL\ep^2)\chi
\\&=O(\ga)\chi+O(\ep)+O(\ga^2\gL \ep^2)\chi +O(\ga\gL\ep^2)\chi
\\&=O(\ga+\ep+\ep^2\ga^2\gL+\ep^2\ga\gL)\chi
\\&=O(\ga+\ep^2\ga^2\gL)\chi,
%\\D_x(\psi-\rho)&=O(\ep^2\ga^2\gL),
\\
\ga\gL\chi (y-\ep h)&=O^+(\ga\ep\gL)\chi,
\\
\gamma^+\cdot D_x\psi+D_y\psi-\beta^+
&=%\gamma_0\cdot D_x\psi +
O(\ep)D_x\psi+
\beta_0 +2\ga\gL\chi (y-\ep h)-\beta^+
\\&=
O(\ep\ga +\ep^3\ga^2\gL)\chi
+O^+(\ep\ga\gL)\chi+O(\ep)
\\&
=O(\ep\ga+\ep^3\ga^2\gL+\ep)\chi+O^+(\ep\ga\gL)\chi
\\&
=O(\ep^2\ga^2\gL+\ep\ga)\chi+O^+(\ep\ga\gL)\chi
\\&
=(O(\ep\ga)+O((\ep\ga)^2\gL)+O^+(\ep\ga\gL)) \chi.
\eald\]
We write the last term above as
\[
(d_1+d_2(\gL)+d_3(\gL))\chi,
\]
where 
\[
d_1=O(\ep\ga),\qquad d_2(\gL)=O((\ep\ga)^2\gL),\qquad d_3(\gL)=O^+(\ep\ga\gL). 
\]
We choose $\gL_0>1$ so that
\[
d_1+d_3=O^+(\ep\ga\gL_0)=:d_4,
\] 
and then choose $\gd_0\in (0,1)$ so that if $\ep\ga\leq\gd_0$, then
\[
d_2(\gL_0)+d_4=O^+(\ep\ga). 
\]
This shows that for $\gL=\gL_0$, if $\ep\ga\leq\gd_0$, then we have
\beq\label{strict-top}
\gamma^+\cdot D_x\psi+\psi_y-\beta^+>0 \ \ \ON \pl_\T\gO_\ep.
\eeq

Henceforth, we fix $\gL=\gL_0$. For brevity of notation, we write $\gL$ to denote $\gL_0$. 

Next, we seek for a kind of concavity of the function $\psi$.  
For this, recalling \eqref{psi_x} and \eqref{psi_y}, we compute the Hessian matrix of $\psi$ in $\gO_\ep$ as follows: let $i,j\in\{1,\ldots,N\}$ and compute that 
\[\bald
%\psi_{x_i}&=-\chi_{x_i}+\beta_{0,x_i}y+\ga\gL\chi_{x_i}(y-\ep h)^2-2\ep\ga\gL\chi (y-\ep h)h_{x_i},
%\\
%\psi_y&=-\beta_0 +\ga\gL\chi (y-\ep h),
%\\
\psi_{x_ix_j}&=-\chi_{x_ix_j}+\beta_{0,x_ix_j} y
+\ga\gL \chi_{x_ix_j}(y-\ep h)^2
-2\ep \ga\gL\chi_{x_i}(y-\ep h)h_{x_j}
\\&\quad -2\ep\ga\gL\chi_{x_j}(y-\ep h)h_{x_i} 
-2\ep\ga\gL\chi (y-\ep h)h_{x_ix_j}
+2\ep^2\ga\gL\chi h_{x_i}h_{x_j},
\\
\psi_{x_iy}&=\beta_{0,x_i}+2\ga\gL \chi_{x_i}(y-\ep h)-2\ep\ga\gL\chi h_{x_i},
\\
\psi_{yy}&=2\ga\gL\chi, 
%\\
%\chi_{x_i}&=\ga \chi s_{x_i}=O(\ga)\chi,
%\\
%\chi_{x_ix_j}&=
%\ga^2\chi s_{x_i}s_{x_j}+\ga\chi s_{x_ix_j}
%\\&=\ga^2\chi s_{x_i}s_{xj}+O(\ga)\chi.
\eald
\] 
Furthermore, noting that $\ep\ga<1, \,\ga>1$ by \eqref{eq9},
$\chi_{x_i}=O(\ga)\chi$, $\chi_{x_ix_j}=O(\ga^2)\chi$, 
and $y=O(\ep)$, we deduce that 
\[\bald 
\psi_{x_ix_j}&=-\ga^2\chi s_{x_i}s_{x_j}+O(\ga)\chi, 
\\
\psi_{yx_i}&=O(1+\ep\ga^2+\ep\ga)\chi =O(\ga)\chi. 
\eald
\]
Thus, we have
\[
D^2\psi=-\ga^2\bmat Ds\otimes Ds &0\\ 0& 0  \emat\chi
+O(\ga)\chi
=-\ga^2 (Ds,0)^\sT(Ds,0)\chi +O(\ga)\chi. 
\]
Combining this and \eqref{positive**} yields 
\[
-\tr A_{\gl\mu} D^2\psi =\chi (\ga^2\tr (Ds,0)A_{\gl\mu}(Ds,0)^\sT+O(\ga))
\geq \chi(\ga^2+O(\ga)).
\]

Note that 
\[
b_{\gl\mu}\cdot D\psi=O(1)|D\psi|=O(\ga)\chi.
\]
Hence, 
\[
-\tr A_{\gl\mu}D^2\psi-b_{\gl\mu}\cdot D\psi-f_{\gl\mu}
\geq \chi(\ga^2+O(\ga)).
\]
We can select $\ga_0>1$ so that for $\ga\geq \ga_0$, 
\[
-\tr A_{\gl\mu} D^2\psi
-b_{\gl\mu}\cdot D\psi  
-f_{\gl\mu} >0. 
%%%\geq \fr{\ga^2}{2}\chi. 
\]

Henceforth, we fix $\ga=\ga_0$. The constant $C_D$ is free to disposal at this point. Since 
$\psi\geq C_D+\beta_0(x)y$, we may fix $C_D>0$ so that \ $\psi>0$\ on 
$\,\ol{\gO_\ep}$ and
\beq\label{alpah-zero}
-\tr A_{\gl\mu} D^2\psi
-b_{\gl\mu}\cdot D\psi+c_{\gl\mu}\psi 
-f_{\gl\mu}
>0. 
\eeq

Now that $\gL=\gL_0>1$ and $\ga=\ga_0>1$ are fixed, the functions $\psi=\ol{\psi_\ep}$ on 
$\ol\gO\tim[-r,r]$ depends only on the parameter $\ep$. Choosing $\ep_1>0$ so that
$\ep_1<\min\{\ep_0, \gd_0/\ga_0\}$, we deduce that for $\ep\in (0,\ep_1)$,  $\ol{\psi_\ep}$ satisfies \eqref{m.1}, \eqref{m.3}, and also \eqref{m.7} for some constant $C$. With similar considerations, we conclude that,  
after redefining $\gL=\gL_0$, $\ga=\ga_0$, 
and $\ep_1$ if necessary,  the pair of functions $\ol{\psi_\ep}$ and $\ul{\psi_\ep}$ satisfies \eqref{m.1}--\eqref{m.7} for $\ep\in(0,\ep_1)$. This ends Step 1. \qed

\begin{remark}\label{after Step1} For the discussion that follows, it will be helpful to note the following: If we review the calculations in Step 1 above, we can see that all the necessary conditions are: 
\begin{enumerate}[(i)]
\item The uniform boundedness of functions $\gs_{\gl\mu}, b_{\gl\mu}, c_{\gl\mu}, f_{\gl\mu}$ (i.e. the first inequality in \eqref{h.2}); 
\item $c_{\gl\mu}$ being non-negative (i.e. \eqref{non-negative}); 
\item the ellipticity condition \eqref{positive} for $A_{\gl\mu}=\gs_{\gl\mu}^\sT\gs_{\gl\mu}$; 
\item the $C^2$-regularity of $\gamma_0, \beta_0$ (guaranteed by \eqref{h.5}) with $\gamma_0$ being zero (i.e. \eqref{gamma00}); 
\item the asymptotic condition:  $\gamma^\pm(x,y)=\pm\gamma_0(x)+O(|y|)$ and $\gb^\pm(x,y)=\pm\gb_0(x)+O(|y|)$; 
\item the inequality $g^-<g^+$ on $\ol\gO$ (i.e. the second condition of \eqref{h.8}).  
\end{enumerate}
\end{remark}

\subsection*{\quad Step 2} Our approach in order to deal with the general situation is to distort the coordinate system a little so that the resulting behavior of $\gamma^\pm$ are the same as those in Step 1.  Henceforth, we assume that functions $\gamma^\pm, \beta^\pm$ (resp., $\gamma_0,\beta_0, k^\pm, l^\pm$)
are defined in $\R^{N}\tim[-1,1]$ (resp., $\R^N$) and satisfy the conditions 
\eqref{cont}, \eqref{h.4}, \eqref{h.5}, \eqref{h.6}, and \eqref{h.7}, with $\ol\gO$ replaced by $\R^N$.  
We may further assume that $\gamma_0$ is bounded in $\R^N$, together with its 
derivatives up to the second order. In this regard, we remind that, by our notation, $f\in C^2(\ol\gO\tim[-1,1])$ means that $f$ are the restrictions to $\ol\gO\tim[-1,1]$ of $C^2$-functions $f$ defined in a neighborhood of $\ol\gO\tim[-1,1]$. 

We may also assume that the functions $g^\pm$ are defined in $\R^N$ as bounded continuous 
functions in $\R^N$ and satisfy $\, g^-(x)<g^+(x)$ for $x\in\R^N$.
The new coordinate system, to be introduced, can be described as the following change of unknowns: 
Let $\gamma:=\gamma_0$ for the notational simplicity. When $u^\ep$ is a solution of \dirichlet, we consider the function $v^\ep$ defined by 
\beq\label{s2.1}
v^\ep(z,y)=u^\ep(z+y\gamma(z),y),
\eeq
which satisfies basically the same boundary condition as in Step 1 on the top and bottom portions of 
boundary of the domain of $v^\ep$. 

Let $0<r<1$ be a small constant. The condition on its size will be revealed later on. 
Consider the function $P(z,y)=(z+y\gamma(z),y)$ in $\R^N\tim[-r,r]$. 

\begin{lemma} \label{s2.lem1}For $r\in(0,1)$ sufficiently small, the mapping $P$ is a $C^2$-diffeomorphism
on $\R^N\tim[-r,r]$. Furthermore, the derivatives, up to second order, of the inverse mapping $P^{-1}$ are bounded 
in $\R^N\tim[-r,r]$. 
\end{lemma} 

\bproof  Since $\gamma\in C^2(\R^N)$ and the Jacobian of the mapping $P$ 
does not vanish on $X:=\R^N\tim[-r,r]$ by assuming $r>0$ sufficiently small, 
we deduce by the Inverse Function Theorem that $P$ is a local $C^2$-diffeomorphsm. 
We may assume that the Jacobian of $P$ is away from a neighborhood of the origin. 
It is obvious that $P$ maps $X$ into $X$. 
We show that, if $r>0$ is sufficiently small, $P$ is a homeomorphism. 
For this, it is enough to show that, if $r>0$ is sufficiently small, then  $P$ is bijective on $
X$ i.e. that  for any $(x,y)\in X$, there exists a unique $z\in\R^N$ such that $x=z+y\gamma(z)$. This is an 
easy consequence of the  Contraction Mapping Theorem, applied to the mapping: 
$z\mapsto x-y\gamma(z)$ in $\R^N$. 

Since the derivatives of $P$ are bounded up to second order and its Jacobian
is away from a neighborhood of the origin, the Inverse Function Theorem assures that 
the derivatives of $P^{-1}$, up to the second order, are bounded. 
\eproof 

Hereafter, we fix $0<r<1$ so that the mapping $P$ on $\R^N\tim[-r,r]$ is a $C^2$-diffeomorphism
and the derivatives $P^{-1}$ up to second order are bounded in $\R^N\tim[-r,r]$.   
It is obvious that $P^{-1}(x,y)=(\gz(x,y),y)$ for all $(x,y)\in\R^N\tim[-r,r]$ and some function 
$\gz\in C^2(\R^N\tim[-r,r])$. 
%%%{\color{magenta}.} 
We note that the derivatives of $\gz$ up to the second order are bounded in $\R^N\tim[-r,r]$. Note also that $z=\gz(x,y)$ if and only if $x=z+y\gamma(z)$. 

We are concerned with the correspondence between ``top boundaries'' for $u^\ep$ and $v^\ep$,
where the relation of $u^\ep$ and $v^\ep$ is given by \eqref{s2.1}. 
We may assume by selecting $\ep_0\in(0,r]$ small enough that $\ep_0 |g^\pm(x)|\leq r$ for $x\in\R^N$. 
Let $\ep\in(0,\ep_0]$. If $u^\ep$ is a solution of \dirichlet, then the ``top boundary'' for $u^\ep$ is given as the graph $y=\ep g^+(x)$, $x\in\ol\gO$.  We remove this restriction, $x\in\ol\gO$, for simplicity, and we instead consider the graph $\gG_\ep^+$ of the function $\ep g^+$ on $\R^N$. 
If $v^\ep$ is given by \eqref{s2.1}, then the ``top boundary'' for $u^\ep$ is mapped to the ``top boundary'', $\Sigma_\ep^+:=P^{-1}(\gG_\ep^+)$, for $v^\ep$. Observe that 
\[\bald
\Sigma_\ep^+&=\{(z,y)\in\R^N\tim[-r,r] \mid P(z,y)\in\gG_\ep^+\}
\\&=\{(z,y)\mid z\in\R^N, y=\ep g^+(z+y\gamma(z))\}. 
\eald 
\]

\begin{lemma}\label{s2.lem2} For $r>0$ sufficiently small, if $0<\ep\leq\ep_0$, then $\Sigma_\ep^+$ is the 
graph of a function $g_\ep^+\in C^2(\R^N)$, with the property 
\[
|g^+_\ep(z)-\ep g^+(z)|\leq C\ep^2\ \ \ \FORALL z\in\R^N,
\]
where $C>0$ is a constant independent of $z, \ep$. Furthermore, for any $(x,y)\in\R^N\tim[-r,r]$, 
if $(z,y)=P^{-1}(x,y)$, then we have 
\[
y>\ep g^+(x) \ \ (\text{resp.,}\ y<\ep g^+(x)) \ \iff \ y>g_\ep^+(z) \ \ (\text{resp.,}\ y<g^+_\ep(z)). 
\]
\end{lemma} 

\bproof For any $z\in\R^N, \ep\in(0,\ep_0]$, consider the $C^1$-function $f(y)=f(y;z,\ep):=y-\ep g^+(z+y\gamma(z))$ on $\R$.  By the choice of $\ep_0$, it follows that $f(r)\geq 0\geq f(-r)$. 
Note that 
\[
f'(y)=1+\ep Dg^+(z+y\gamma(z))\cdot \gamma(z)\geq 1-r\|Dg^+\|_\infty \|\gamma\|_\infty. 
\]
The right-hand side of the above can be assumed to be larger than $1/2$. 
By the Implicit Function Theorem, there exists a $C^1$-function $g^+_\ep$ 
such that $g^+_\ep(z)\in[-r,r]$ and  
for $(z,y)\in \R^N\tim[-r,r]$, 
\[
y=g^+_\ep(z) \ \iff \ f(y;z,\ep)=0. 
\] 
Obviously, the function $f(y)$ is strictly increasing in $\R$ and the
above condition implies that 
\[\bald
y>g^+_\ep(z) \ \ (\text{resp., }y<g^+_\ep(z)) \ &\iff\ f(y;z,\ep) >0 \ \ (\text{resp., } f(y;z,\ep)<0)
\\
& \iff\ y>\ep g^+(z+y\gamma(z)) \ \ (\text{ resp., } y<\ep g^+(z+y\gamma(z))). 
\eald
\]
For any $z\in\R^N$, if we write $y=g_\ep^+(z)$, then 
\[
y=\ep g^+(z+y\gamma(z)), \qquad
|y|\leq \ep\|g^+\|_\infty,
\]
and, moreover, 
\begin{align}\notag
|\ep g^+(z)-g^+_\ep(z)|&=\ep|g^+(z)-g^+(z+y\gamma(z))|
\\&\leq \ep \|Dg^+\|_\infty\|\gamma\|_\infty |y|
\leq \ep^2\|Dg^+\|_\infty\|g^+\|_\infty \|\gamma\|_\infty. \qedhere 
\end{align}

Similalrly, we denote by $\gG_\ep^-$ the graph of the function $\ep g^-$ and set 
\[
\Sigma_\ep^-:=P^{-1}(\gG_\ep^-)=\{(z,y)\mid z\in\R^N,y=\ep g^-(z+y\gamma(z))\}. 
\]
We then have the following lemma, whose proof can be done in the same way as that of Lemma 
\ref{s2.lem2} and is left to the reader.

\begin{lemma}\label{s2.lem3} For $r>0$ sufficiently small, if $0<\ep\leq\ep_0$, then $\Sigma_\ep^-$ is the 
graph of a function $g_\ep^-\in C^2(\R^N)$, satisfying the inequality 
\[
|g^-_\ep(z)-\ep g^-(z)|\leq C\ep^2\ \ \ \FORALL z\in\R^N,
\]
where $C>0$ is a constant independent of $z, \ep$. Furthermore, for any $(x,y)\in\R^N\tim[-r,r]$, 
if $(z,y)=P^{-1}(x,y)$, then we have 
\[
y>\ep g^-(x) \ \ (\text{resp.,} \ y<\ep g^-(x)) \ \iff \ y>g_\ep^-(z) \ \ (\text{resp.,} \ y<g^-_\ep(z)). 
\]
\end{lemma}

\begin{lemma}\label{s2.lem4} Let $0<r<1$ be sufficiently small so that $P$ is a $C^2$-diffeomorphism on $\R^N\tim[-r,r]$.Then, for any $(x,y)\in\R^N\tim[-r,r]$, 
\[
|P^{-1}(x,y)-(x,y)|\leq r\|\gamma\|_\infty. 
\]
\end{lemma}

\bproof Observe that
\begin{align*}
\sup_{(x,y)\in\R^N\tim[-r,r]}|P^{-1}(x,y)-(x,y)|
&=\sup_{(z,y)\in\R^N\tim[-r,r]} |(z,y)-P(z,y)|
\\&=\sup_{(z,y)\in\R^N\tim[-r,r]}|y||\gamma(z)|\leq r\|\gamma\|_\infty.\qedhere
\end{align*}
\eproof

\eproof

\subsection*{\quad Step 3} This is the final step of the proof of Lemma \ref{main.lem}. 
Let $r\in(0,1)$ be a constant sufficiently small, so that Lemmas \ref{s2.lem1}--\ref{s2.lem4} 
are valid. Let $P$ be the $C^2$-diffeomorphism on $\R^N\tim[-r,r]$. 
For notational convenience, we also write $Q$ for $P^{-1}$. We remark that the mappings 
$z\mapsto P(z,0)$ and $x\mapsto Q(x,0)$ are the identity mapping in $\R^N$. 
 Let $g_\ep^\pm$ be $C^1$-functions in $\R^N$, from Lemmas \ref{s2.lem2} and \ref{s2.lem3}. We always assume that 
$\ep_0\|g^\pm\|_\infty\leq r$, which assures that 
\[
\ol{\gO_\ep}\subset \{(x,y)\in\R^{N+1}\mid \ep g^-(x)\leq y\leq \ep g^+(x)\}\subset \R^N\tim[-r,r]. 
\]
Here and hereafter, we assume that $0<\ep\leq\ep_0$. It follows that 
$Q(\ol{\gO_\ep})\subset \R^N\tim[-r,r]$. By Lemmas \ref{s2.lem2} and \ref{s2.lem3}, 
we see that 
\[\bald
Q(\gO_\ep)&\subset \{(z,y)\in\R^{N+1}\mid g_\ep^-(z)<y<g_\ep^+(z)\},\qquad
\\ Q(\ol{\gO_\ep})&\subset \{(z,y)\in\R^{N+1}\mid g_\ep^-(z)\leq y\leq g_\ep^+(z)\}.
\eald
\]
Since $Q(\gO_\ep)$ is, in general, not a truncated cylinder, we wish to replace it with a 
truncated cylinder 
\[
\Hat\gO_{\ep}=\{(z,y)\in\Hat\gO\tim\R\mid g_\ep^-(z)<y<g_\ep^+(z)\},
\]
where $\Hat\gO$ is a bounded $C^2$-domain and a neighborhood of $\ol\gO$.  
We extend the domain of the functions $\gs_{\gl\mu}, b_{\gl\mu}, c_{\gl\mu},$
and $f_{\gl\mu}$ so that \eqref{h.2} and \eqref{non-negative} are satisfied with 
$\R^N$ in place of $\gO$ and with newly selected  
$C_F, \go_F$. 
Let $s$ be the function from \eqref{positive}, and we extend its domain to 
$\R^N$ so that $s\in C^2(\R^N)$. 
  
We may assume by choosing $\Hat\gO$ small enough that  
%\[(Ds(x),0)A_{\gl\mu}(x,0)\not=0 %(Ds(x),0)^\sT>0 \ \ \FOR x\in\ol{\widetilde \gO},
%\]
%and we may furthermore assume by choosing $r>0$ smaller if necessary that   
\beq \label{positive+}
\inf_{(\gl,\mu)\in L\tim M}\left|(Ds(x),-Ds(x)\gamma^\sT(x))A_{\gl\mu}(x,0)\right|>0 %(Ds(x),0)^\sT>0 
\ \ \FOR x\in\ol{\Hat \gO}.  %\tim[-r,r].
\eeq 

Thanks to Lemma \ref{s2.lem4}, we may assume by repalcing $r>0$ with a smaller $r>0$ if necessary   
that $Q(\ol{\gO_\ep})\subset \Hat\gO\tim\R$ for every $\ep\in(0,\ep_0]$. Thus, we deduce that 
\[Q(\ol{\gO_\ep})
\subset \Hat\gO_\ep.
\]
By definition, we have $\pl_\T\Hat \gO_\ep=\{(z,y)\in\ol {\,\Hat\gO\,}\tim\R\mid y=g_\ep^+(z)\}$ 
and $\pl_\B\Hat\gO_\ep=\{(z,y)\in\ol{\,\Hat\gO\,}\tim\R\mid y=g_\ep^-(z)\}$. 
Observe by Lemmas \ref{s2.lem2} and \ref{s2.lem3} that 
\[Q(\pl_\T\gO_\ep)\subset \pl_\T \Hat\gO_\ep \ \ \AND \ \ Q(\pl_\B\gO_\ep)\subset \pl_\B \Hat\gO_\ep. 
\]
We introduce functions $\hat \gamma^\pm$, $\hat\gb^\pm$ on $\R^N\tim[-r,r]$ by setting
\beq\label{hat-gamma-beta}
\hat \gamma^\pm(z,y):=\gamma^\pm(P(z,y))\left(DP(z,y)^\sT\right)^{-1},\qquad
\hat \gb^\pm(z,y):=\gb^\pm(P(z,y)).
\eeq
Define $\Hat F\in C(\cS(N+1)\tim\R^{N+1}\tim\R\tim (\R^N\tim[-r,r]))$ by
\beq\label{hat_F}\Hat F(X,p,t,(z,y)):=F(Y,q,t,P(z,y)),
\eeq  
where $F$ is considered as a function on $\cS(N+1)\tim\R^{N+1}\tim\R\tim (\R^N\tim[-r,r])$ defined by \eqref{h.1} with $\R^N$ in place of $\gO$; 
\beq\label{s3.Yq}
Y=DQ^\sT(P(z,y))XDQ(P(z,y))+p\cdot D^2Q(P(z,y)),\qquad 
q=p DQ(P(z,y)); 
\eeq
and, for $Q(\xi)=(Q_1(\xi),\ldots,Q_{N+1}(\xi))$, $p\cdot D^2Q(\xi)$ denotes the $(N+1)\tim(N+1)$ matrix, with its $(i,j)^{\,\text{th}}$ entry is given by 
\[
\sum_{k=1}^{N+1}p_k Q_{k,\xi_i\xi_j}.
\] 
With the correspondence \eqref{s3.Yq}, we set 
\[
\Hat F_{\gl\mu}(X,p,t,(z,y)):=F_{\gl\mu}(Y,q,t,P(z,y)),
\]
and  we have
\[
\Hat F(X,p,t,(z,y))=\inf_{\gl\in L}\sup_{\mu\in M}\Hat F_{\gl\mu}(X,p,t,(z,y)).
\]

Let $(x,y)\in\R^N\tim[-r,r]$ and set $(z,y)=Q(x,y)$.
 Note that
$(x,y)=P(z,y)$ and 
\[
DP(z,y)=\bmat I_N+y D\gamma(z)&\gamma(z)^\sT\\ 0&1\emat, 
\]
which implies  
\beq\label{s3.15} %%\bald
DQ(x,y)=DP(z,y)^{-1}=\bmat (I_N+yD\gamma)^{-1} & -(I_N+yD\gamma)^{-1}\gamma^\sT \\ 0&1\emat.  
%% \\
%%\left(DP^\sT(z,y)\right)^{-1}&=\left(DP(z,y)^{-1}\right)^\sT
%%=\bmat (I_N+yD\gamma^\sT)^{-1} &0 \\ 
%%%-\gamma (I_N+yD\gamma^\sT)^{-1}& 1\emat.
%\eald
\eeq
For $(z,y)\in \R^N\tim[-r,r]$, we set 
\[
R(z,y)=\bmat (I_N+yD\gamma(z))^{-1} & -(I_N+yD\gamma(z))^{-1}\gamma^\sT(z) \\ 0&1\emat,
\]
and note that $DQ\circ P=R$.

Let $(z,y)\in \R^N\tim[-r,r]$ and $(X,p,t)\in\cS(N+1)\tim\R^{N+1}\tim\R$, and let 
$(Y,q)\in\cS(N+1)\tim\R^{N+1}$ be given by \eqref{s3.Yq}. Set $(x,y)=P(z,y)$. It follow from the above computation 
that 
\[\bald
\tr A_{\gl\mu}(P(z,y))Y&=\tr A_{\gl\mu}(x,y)R^\sT(z,y)X R(z,y)+\tr A_{\gl\mu}(x,y) p\cdot D^2Q(x,y)
\\&= \tr R(z,y)A_{\gl\mu}(x,y)R^\sT(z,y)X +\tr A_{\gl\mu}(x,y) p\cdot D^2Q(x,y),
\eald
\] 
\[
R(z,y)A_{\gl\mu}(x,y)R^\sT(z,y)=\left(\gs_{\gl\mu}(x,y)R^\sT(z,y)\right)^\sT \gs_{\gl\mu}(x,y)R^\sT(z,y),
\]
\[
\tr A_{\gl\mu}(x,y) p\cdot D^2Q(x,y)=\sum_{i=1}^{N+1} p_i \tr A_{\gl\mu}(x,y)D^2 Q_i(x,y),
\]
and also,
\[\bald
b_{\gl\mu}(P(z,y))\cdot q&=b_{\gl\mu}(x,y)\cdot pDQ(P(z,y))=b_{\gl\mu}(x,y)(pR(z,y))^\sT
\\&=b_{\gl\mu}(x,y)R^\sT(z,y) p^\sT=b_{\gl\mu}(x,y)R^\sT(z,y)\cdot p. 
\eald
\]
We set 
\beq\label{s3.18}\left\{\,\bald
d_{\gl\mu}:&=\left(\tr A_{\gl\mu}D^2Q_1,\ldots, \tr A_{\gl\mu}D^2Q_{N+1}\right), 
\\
\hat\gs_{\gl\mu}:&= \gs_{\gl\mu}\circ P \,R^\sT,
\\
\hat b_{\gl\mu}:&= b_{\gl\mu}\circ P +d_{\gl\mu}\circ P,
\\
\hat c_{\gl\mu}:&=c_{\gl\mu}\circ P, 
\\
\hat f_{\gl\mu}:&=f_{\gl\mu}\circ P,
\eald\right.
\eeq
to obtain 
\begin{align}\label{s3.19}
\Hat F_{\gl\mu}(X,p,t,(z,y))&=-\tr (\hat\gs_{\gl\mu}^\sT \hat \gs_{\gl\mu})(z,y) X
-\hat b_{\gl\mu}(z,y)\cdot p
%\\&\quad 
+\hat c_{\gl\mu}(z,y) t-\hat f_{\gl\mu}(z,y) 
\end{align}
for $(X,p,t,(z,y))\in \cS(N+1)\tim \R^{N+1}\tim\R\tim \R^N\tim[-r,r]$.
Observe that the functions $\hat\gs_{\gl\mu}, \hat b_{\gl\mu}, \hat c_{\gl\mu}, \hat f_{\gl\mu}$ 
are bounded in $\R^N\tim[-r,r]$ uniformly for $(\gl,\mu)\in L\tim M$ and that $\hat\gs_{\gl\mu}\geq 0$.   We also set 
$\Hat A_{\gl\mu}:=\hat\gs_{\gl\mu}^\sT\hat\gs_{\gl\mu}$.

Next, we examine some properties of $\hat \gamma^\pm$ and $\hat\beta^\pm$. 
For this, we compute that for $(z,y)=Q(x,y)\in\ol{\Hat \gO} \tim[-r,r]$, 
\[%%%\label{s3.16}
\bald
\hat \gamma^\pm (z,y)&=(\gamma_1^\pm(P(z,y)),\pm 1)\bmat \left(I_N+y D\gamma^\sT\right)^{-1} 
&0\\ -\gamma (I_N+yD\gamma^\sT)^{-1} &1\emat
\\&=\left(\left(\gamma_1^\pm(P(z,y))\mp\gamma(z)\right)\left(I_N+yD\gamma^\sT(z)\right)^{-1}, \pm 1\right),
\eald
\]
and
\[\bald
\left(\gamma_1^\pm(P(z,y))\mp\gamma(z)\right)\left(I_N+yD\gamma^\sT(z)\right)^{-1}
&
=(\gamma_1^\pm(z+y\gamma(z),y)\mp\gamma(z))(I_N+O(|y|)
\\&
=\pm\gamma(z+y\gamma(z))\mp\gamma(z))+O(|y|)
=O(|y|).
\eald
\]
Hence, we find that if we set $\hat\gamma_0=0$ in $\ol{\Hat\gO}$, then 
\beq\label{s3.17}
\hat \gamma^\pm(z,y)=(\hat\gamma_0(z)+O(|y|),\pm 1) \ \ \FOR (z,y)\in \ol{\Hat\gO}\tim[-r,r]. 
\eeq
Similarly, we deduce that for $(z,y)\in \ol{\Hat\gO}\tim[-r,r]$, 
\beq\label{hat_beta}
\hat\beta^\pm(z,y)=\beta^\pm(z+y\gamma(z),y)=\pm \beta_0(z)+O(|y|). 
\eeq

Observe that
for $z\in \ol{\Hat\gO}$,
\[\bald
(Ds(z),-Ds(z)\hat\gamma_0^\sT(z))\hat \gs_{\gl\mu}^\sT(z,0)
&=(Ds(z),0)R(z,0)\gs_{\gl\mu}^\sT(z,0)
\\&=(Ds(z),0)\bmat I_N &-\gamma^\sT(z)\\ 0&1\emat \gs_{\gl\mu}^\sT(z,0)
\\&=\left(Ds(z), Ds(z)\gamma^\sT(z)\right) \gs_{\gl\mu}^\sT(z,0),
\eald
\]
from which we obtain
\[\bald
&(Ds(z),-Ds(z)\hat\gamma_0^\sT(z))\Hat A_{\gl\mu}^\sT(z,0)
(Ds(z),-Ds(z)\hat\gamma_0^\sT(z))^\sT 
\\&\quad =\left(Ds(z), Ds(z)\gamma^\sT(z)\right) A_{\gl\mu}(z,0)
\left(Ds(z), Ds(z)\gamma^\sT(z)\right)^\sT.
\eald
\]
and, by \eqref{positive+}, 
\beq\label{hat-positive}
\inf_{(\gl,\mu)\in L\tim M}
\left|(Ds(z),-Ds(z)\hat\gamma_0^\sT(z)){\Hat A_{\gl\mu}}^\sT(z,0)\right|>0
\ \ \FOR z\in \ol{\Hat\gO}. 
\eeq

Now, in view of Remark \ref{after Step1}, we are in a good position to apply the argument of Step 1 
for the collection of the functions $\Hat F, \hat\gamma^\pm, \hat\beta^\pm$ and the domain $\Hat\gO$. By \eqref{hat-positive}, \eqref{s3.17}, \eqref{hat_beta}, this collection 
has the properties (iii), (iv), (v) stated in Remark \ref{after Step1}. Also, the properties (i) and (ii)  are 
already checked. The Step 1 argument will be done with the thin domains $\Hat\gO_\ep$,
and we have to replace $\ep g^\pm$ with $g_\ep^\pm$. This does not create any difficulty 
and we need to notice that for any function $h\in C^2(\R^N)$ satisfying 
$g^-<h<g^+$ in $\R^N$,  if $\ep>0$ is sufficiently small, then 
\[
\min\{g_\ep^+-\ep h,\ep h-g_\ep^-\}\geq \ep \gd \ \ \ON \ol{\Hat\gO}
\] 
for some constant $\gd>0$ (see \eqref{s2.lem3}). With regard to $g_\ep^\pm$, this 
observation is enough to carry out Step 1.   

Thus, the argument in Step 1 guarantees the existence of  functions $\ol{w_\ep},\,\ul{w_\ep}\in 
C^2(\ol{\Hat\gO}\tim[-r,r])$, with $\ep\in(0,\ep_1)$, where $\ep_1$ is a constant in $(0,\ep_0)$ , such that for  any $\ep\in(0,\ep_1)$,  
\begin{align}
\label{s3.1}&\hat \gamma^+\cdot D\ol{w_\ep}-\tilde \gb^+>0 \ \ \ON \pl_\T \Hat\gO_\ep,
\\\label{s3.2}&\hat \gamma^-\cdot D\ol{w_\ep}-\tilde \gb^->0 \ \ \ON \pl_\B\Hat\gO_\ep,
\\\label{s3.3}&\Hat F(D^2\ol{w_\ep},D\ol{w_\ep},\ol{w_\ep},(z,y))>0 \ \ \ON \ol{\Hat \gO_\ep},
\\\label{s3.4}&\hat \gamma^+\cdot D\ul{w_\ep}-\tilde \gb^+<0 \ \ \ON \pl_\T\Hat\gO_\ep,
\\\label{s3.5}&\hat \gamma^-\cdot D\ul{w_\ep}-\tilde \gb^-<0 \ \ \ON \pl_\B\Hat\gO_\ep,
\\\label{s3.6}&\Hat F(D^2\ul{w_\ep},D\ul{w_\ep},\ul{w_\ep},(z,y))<0 \ \ \ON \ol{\Hat \gO_\ep},
\\\label{s3.7}&C>%{\color{magenta} 
\ol{w_\ep}  %}
>\ul{w_\ep}>-C \ \ \ \ON \ol{\,\Hat\gO\,}\tim[-r,r], 
\end{align}
for some constant $C>0$. Noting that $Q(\ol{\gO_\ep})
\subset \Hat\gO_\ep\subset \Hat\gO\tim[-r,r]$, we define $\ol{\psi_\ep}, \ul{\psi_\ep}\in C^2(\ol{\gO_\ep})$ 
by $\ol{\psi_\ep}:=\ol{w_\ep}\circ Q$ and $\ul{\psi_\ep}:=\ul{w_\ep}\circ Q$. Then, from 
\eqref{s3.1}--\eqref{s3.7}, we obatin 
\begin{align}
\label{s3.8}&\gamma^+\cdot D\ol{\psi_\ep}-\gb^+>0 \ \ \ON \pl_\T\gO_\ep,
\\\label{s3.9}&\gamma^-\cdot D\ol{\psi_\ep}-\gb^->0 \ \ \ON \pl_\B\gO_\ep,
\\\label{s3.10}&F(D^2\ol{\psi_\ep},D\ol{\psi_\ep},\ol{\psi_\ep},(x,y))>0 \ \ \ON \ol{\gO_\ep},
\\\label{s3.11}&\gamma^+\cdot D\ul{\psi_\ep}-\gb^+<0 \ \ \ON \pl_\T\gO_\ep,
\\\label{s3.12}& \gamma^-\cdot D\ul{\psi_\ep}-\gb^-<0 \ \ \ON \pl_\B\tilde\gO_\ep,
\\\label{s3.13}&F(D^2\ul{\psi_\ep},D\ul{\psi_\ep},\ul{\psi_\ep},(z,y))<0 \ \ \ON \ol{\gO_\ep},
\\\label{s3.14}&C> %{\color{magenta}
\ol{\psi_\ep}  %}
>\ul{\psi_\ep}>-C \ \ \ \ON \ol{\gO_\ep},
\end{align}
To check these claims, let $(x,y)\in\ol{\gO_\ep}$ and set $(z,y)=Q(x,y)$.
Assume that 
$(x,y)\in \pl_\T\gO_\ep$. Note that $(z,y)\in\pl_\T\Hat\gO_\ep$ and from \eqref{s3.1} that 
$\hat\gamma^+\cdot D\ol{w_\ep}-\hat \gb^+>0$ at $(z,y)$. Observe that
\[\bald
D\ol{\psi_\ep}(x,y)
&=D\ol{w_\ep}\circ Q(x,y)
=D\ol{w_\ep}(z,y)DQ(x,y), 
\\
\gamma^+(x,y)\cdot D\ol{\psi_\ep}(x,y)-\gb^+(x,y)
&=\gamma^+(P(z,y))(D\ol{\psi_\ep}(x,y))^\sT-\gb^+(P(z,y))
\\&=\gamma^+(P(z,y))(D\ol{w_\ep}(z,y)DQ(x,y))^\sT-\gb^+(P(z,y))
\\&=\gamma^+(P(z,y))(D\ol{w_\ep}(z,y)DP(z,y)^{-1})^\sT-\hat \gb^+(z,y)
\\&=\hat \gamma^+(z,y)(D\ol{w_\ep}(z,y)^\sT-\hat \gb^+(z,y),
\eald
\] 
and, therefore, $\gamma^+(x,y)\cdot D\ol{\psi_\ep}(x,y)-\gb^+(x,y)>0$, proving \eqref{s3.8}.  
Similarly, we find that \eqref{s3.9}, \eqref{s3.11}, \eqref{s3.12} are all satisfied. Now, assume
$(x,y)\in\ol{\gO_\ep}$, which implies that $(z,y)\in\ol{\Hat\gO_\ep}$
Since $\ol{\psi_\ep}(x,y)=\ol{w_\ep}(Q(x,y))$, we have 
\[\bald
D^2\ol{\psi_\ep}(x,y)&=DQ^\sT(x,y) D^2w_\ep(z,y) DQ(x,y)+D\ol{w_\ep}(z,y)\cdot D^2Q(x,y)
\\&=DQ(P(z,y))^\sT D^2w_\ep(z,y)DQ(P(z,y))+D\ol{w_\ep}(z,y)\cdot D^2Q(P(z,y)),
\eald
\]
and 
\[\bald
&F(D^2\ol{\psi_\ep}(x,y),D\ol{\psi_\ep}(x,y), \ol{\psi_\ep}(x,y),(x,y))
\\&=\Hat F(D^2\ol{w_\ep}(z,y),D\ol{w_\ep}(z,y),\ol{w_\ep}(z,y),(z,y))>0,
\eald
\]
which shows that \eqref{s3.10} is valid. A parallel computation assures that \eqref{s3.13} holds. 
Inequality \eqref{s3.14} follows immediately from \eqref{s3.7}. 

The proof of Lemma \ref{main.lem} is complete. \qed

\section{Proof of the main theorem}

\subsection{Comparison principle for \dirichletzero}

Let $G_{\gl\mu}^0$ denote the homogeneous part of $G_{\gl\mu}$. That is, for 
$(X,p,r,x)\in\cS(N)\tim\R^N\tim\R\tim\ol\gO$, 
\[
G_{\gl\mu}^0(X,p,r,x)=-\tr \Tilde A_{\gl\mu}(x)X-\tilde b_{\gl\mu}(x)\cdot p 
+\tilde c_{\gl\mu}(x)r.
\]
Note that for $(X_i,p_i,r_i,x)\in\cS(N)\tim\R^N\tim\R\tim\ol\gO$, with $i=1,2$,
\beq\label{additive}
G(X_1,p_1,r_1,x)-G(X_2,p_2,r_2,x)
\leq \sup_{(\gl,\mu)\in L\tim M}G^0_{\gl\mu}(X_1-X_2,p_1-p_2,r_1-r_2,x). 
\eeq

\begin{lemma} \label{perturbation}Under the hypotheses of Proposition \ref{main.prop},
there is a function $\psi\in C^2(\ol\gO)$ such that 
\[
G_{\gl\mu}^0(D^2\psi,D\psi,0,x)\leq -1 \ \ \IN \gO.
\]
\end{lemma}

\bproof Let %%{\color{magenta} 
$s\in C^2(\ol\gO,\R)$ %} 
be from the assumption \eqref{positive}. We may assume that for $(\gl,\mu,x)\in L\tim M\tim\ol\gO$, 
\[
Ds(x)\Tilde A_{\gl\mu}(x)Ds^\sT(x)\geq 1. 
\]
We may assume, by adding a constant to $s$ if necessary, that $s(x)\geq 0$ in $\ol\gO$. 
For $\ga\geq 1$, we set $\psi(x)=e^{\ga s(x)}$ and observe that 
$\psi\geq 1$ in $\ol\gO$, $\psi\in C^2(\ol\gO)$, and 
\[
D\psi=\ga Ds \psi,\quad D^2\psi=\ga^2 \psi Ds\otimes Ds+\ga \psi D^2s. 
\] 
Moreover, using \eqref{l.1}, we compute that
\[\bald
G^0_{\gl\mu}(D^2\psi(x),D\psi(x),0,x)
&=\psi(x)(-\ga^2 Ds(x)\Tilde A_{\gl\mu}(x)Ds^\sT(x)+O(\ga))
\\&\leq \psi(x)(-\ga^2 +O(\ga)).
\eald
\]
By selecting $\ga>1$ large enough, we may assume that the factor, $-\ga^2+O(\ga)$, 
on the right-hand side of the last inequality is less than  $-1$, which yields
that 
\[G^0_{\gl\mu}(D^2\psi(x),D\psi(x),0,x)\leq -1\qquad \FOR\ x\in\ol\gO.
\] 
The function $\psi$ has the required property. 
\eproof

We recall a proposition from \cite{BBI2}. 
\begin{lemma}[{\cite[Lemma 5.3 and Remark 5.4]{BBI2}}] 
\label{pointwise}
 Assume \eqref{h.1}--\eqref{h.3}, \eqref{h.4}--\eqref{h.6}, and \eqref{b-positive}.
Then, the Dirichlet condition in \emph{\dirichletzero} can be understood in the classical sense. 
\end{lemma}

In \cite[Lemma 5.3]{BBI2},  it is assumed instead of \eqref{non-negative} that 
\beq\label{c-positive}
c_{\gl\mu}(x)\geq \ga_0\ \ \ \FORALL \, (x,\gl, % {\color{magenta} , }
\mu)\in \ol\gO\tim L\tim M\ 
\text{ and \ for some constant }\ga_0>0.
\eeq 
However, 
by reviewing the proof of \cite[Lemma 5.3]{BBI2}, 
one easily realizes that \eqref{non-negative} is indeed enough. 

The comparison principle for \dirichletzero is stated as follows. 

\begin{lemma}\label{comparison} Under the hypotheses of Theorem \ref{main.thm}, let $v\in\USC(\ol\gO,\R)$ and $w\in \LSC(\ol\gO,\R)$
be, respectively,  viscosity sub- and super-solution to 
\[
G(D^2u,Du,u,x)=0 \ \ \IN \gO.
\] 
If $v\leq w$ on $\pl\gO$, then $v\leq w$ in $\gO$. 
\end{lemma} 

\bproof[Outline of proof] Let $\psi\in C^2(\ol\gO)$ be the function from Lemma \ref{perturbation}.   
We may assume by adding a constant that $\psi\leq 0$ in $\ol\gO$. Let $\gd>0$ be a constant. 
We set $v_\gd:=v+\gd\psi$ and $w_\gd:=w-\gd\psi$. Using \eqref{additive}, we easily deduce
that $v_\gd$ and $w_\gd$ are sub- and super-solutions to 
$G=-\gd$ and $G=\gd$ in $\gO$, respectively. We also have $v_\gd\leq w_\gd$ on $\pl\gO$. 
 Thanks to \eqref{l.3} and the nonegativity of $\tilde c_{\gl\mu}$, we may apply the standard 
 comparison principle (for instance, \cite{CIL}), to yield $v_\gd\leq w_\gd$ in $\gO$, from which 
 we conclude that $v\leq w$ in $\gO$.
\eproof

\subsection{Proof of the theorem} 
\bproof[Proof of Theorem \ref{main.thm}] Let $S_\ep$ be the set of all viscosity solutions to 
\dirichlet. 
Propostion \ref{main.prop} assures that there are constants $\ep_1\in (0,1)$ and $C_0>0$ such 
that, for every $\ep\in (0, \ep_1)$,  $S_\ep\not=\emptyset$ and $\sup_{\ol{\gO_\ep}}|v|\leq C_0$ for all 
$v\in S_\ep$.   In particular, assertion (i) is valid. For $x\in\ol\gO$, set 
\[\bald
u^+(x):&=\lim_{\ep\to 0^+}\sup \{v(z,y)\mid v\in S_\gd, (z,y)\in \ol{\gO_\gd}, |x-z|<\gd,0<\gd<\ep \},
\\ u^-(x):&=\lim_{\ep\to 0^+}\inf \{v(z,y)\mid v\in S_\gd, (z,y)\in \ol{\gO_\gd}, |x-z|<\gd, 0<\gd<\ep \}.
\eald
\]
The half-limits $u^\pm$ are bounded functions in $\ol\gO$ and $u^+, u^-\in \USC(\ol\gO)$.  According to 
\cite[Theorem 4.1]{BBI2}, whose conclusion remains the same after \eqref{c-positive} is relaxed to \eqref{non-negative}, we find that $u^+$ and $u^-$ are viscosity sub- and super-solutions
to \dirichletzero, respectively. Now, Lemmas \ref{pointwise} and \ref{comparison} assure that 
$u^+(x)\leq\gb(x,0)\leq u^-(x)$ for  $x\in\pl\gO$ and $u^+\leq u^-$ in $\gO$.
Furthermore, we see that $u^+=u^-$ in $\ol\gO$. If we write $u^0:=u^+=u^-$, then 
$u^0\in C(\ol\gO)$, $u^0$ is a viscosity solution to \dirichletzero, and 
\[
\lim_{\ep\to 0^+}\sup_{u\in S_\ep}\sup_{(x,y)\in\ol{\gO_\ep}}|u(x,y)-u^0(x)|=0.
\]     
This shows assertion (iii) is valid. Noting by Lemmas \ref{pointwise} and \ref{comparison} the 
uniqueness of viscosity solution to \dirichletzero, we see that assertion (ii) holds. 
\eproof

\end{document}